# AN UPPER BOUND FOR THE PRIME GAP

Ya-Ping Lu and Shu-Fang Deng

**Abstract**: We showed that the prime gap for a prime number $p$ is less than or equal to the prime count of the prime number, or $g(p) \leq \pi(p)$.

***Definition 1**: The prime gap, $g(p_n)$, is defined as the number of integers greater than the nth prime number $p_n$ and less than or equal to the next prime number $p_{n+1}$:*

$$g(p_n) := \#\{N: p_n < N \leq p_{n+1}\}$$

*Thus, the prime gap equals to the difference between two consecutive prime numbers, $g(p_n) = p_{n+1} - p_n$.*

The smallest prime gap occurs between the first and the second prime numbers, 2 and 3, and $g(p_1) = p_2 - p_1 = 1$. For twin primes $p_n$ and $p_n + 2$, $g(p_n) = 2$. As $n$ increases, larger prime gaps are expected to appear and, as a matter of fact, prime gap can be arbitrarily large. For example, the $m - 1$ consecutive integers $m! + 2, m! + 3, m! + 4, \ldots, m! + m$ are all composites[1]. If $p$ is the largest prime number less than or equal to $m! + 1$, then the prime gap $g(p) \geq m$.

From the prime number theorem, the number of primes less than or equal to $p_n$ is approximately $p_n/ \log p_n$, or $p_n \sim n \log p_n$. So, on average, the prime gap between two prime numbers $p_n$ and $p_{n+1}$ is $\log p_n$.

Cramér[2] showed that if the Riemann Hypothesis holds,

$$g(p) < k p^{1/2} \log p.$$

Bertrand's postulate[3] states that for $n \geq 1$, $p_{n+1} < 2p_n$, or

$$g(p_n) < p_n$$

Bertrand's postulate was proved by Chebyshev[4] and so the postulate is also called the Bertrand–Chebyshev theorem or Chebyshev's theorem.

Based on Ramanujan's work[5] it is proved that $2p_{n-m} > p_n$ for $n > k$, where $k = \pi(p_k) = \pi(R_m)$, and $R_m$ is the $m$-th Ramanujan prime. This means that, for $n > \pi(R_m)$,

$$g(p_n) \leq 2p_{n-m} - p_n$$

One can show, from the prime number theorem, that for every real number $e > 0$ and there is some integer $m_0$ such that there is always a prime $p$ satisfying $m < p < (1 + e)m$ for every $m > m_0$. This shows that, for all $n > n_0$,

$$g(p_n) < ep_n.$$

Some specific pairs of ($e$, $n_0$) are (1/5, 9), (1/13, 118), and (1/16597, 2010760)[6].

Legendre's conjecture states that, for every $m > 1$, there is a prime $p$, such that $m^2 < p < (m+1)^2$. Oppermann's conjecture[7] states that, for every integer $m > 1$, there is at least one prime number between $m(m-1)$ and $m^2$, and at least another prime between $m^2$ and $m(m+1)$. If Oppermann's conjecture is true, there would be at least four prime numbers between $(p_n)^2$ and $(p_{n+1})^2$ for every $n \geq 2$, and the largest possible gaps between two consecutive prime numbers could be, as stated by Andrica's conjecture[8],

$$g(p_n) < 2\sqrt{p_n} + 1$$

which suggests that $g(p_n) = \mathcal{O}(p_n{}^\theta)$ with $\theta < 1/2$ will suffice to prove Andrica's conjecture. However, all values of $\theta$ proved so far are larger than 1/2 and the best unconditional result is $\theta < 21/40$ by R.C. Baker et al.[9].

Table 1 lists some of the prime gaps for $p$ up to 436273009, including the first 30 maximal prime gaps[10], defined as the prime gaps larger than all gaps between smaller primes. We see that the prime gap $g_n$ (Column 3) is smaller than or equal to the prime count of the prime number (Column 1), or $g(p_n) \leq \pi(p_n) = n$.

***Theorem 1***: *The prime gap of a prime number p is less than or equal to the prime count of the prime number, or*

$$g(p) \leq \pi(p).$$

*Proof:* Let $q$ be the next prime number following the prime number $p$. From Definition 1, we have $g(p) = q - p$. Since $\pi(q) - \pi(p) = 1$,

$$\pi(p + g(p)) - \pi(p) = 1.$$

As the prime counting function is non-decreasing, Theorem 1 can be rephrased as "*There is at least one prime number in the range of $(p, p + \pi(p)]$*", or

$$\pi(p + \pi(p)) - \pi(p) \geq 1$$

Dusart[11] showed that the number of prime numbers less than or equal to x is bounded by

$$\pi(x) \leq \frac{x}{\log x - 1.1} \qquad \text{if } x \geq 60184$$

and

$$\pi(x) \geq \frac{x}{\log x - 1} \qquad \text{if } x \geq 5393$$

For $p > 60184$, the number of prime numbers between $p$ and $p + \pi(p)$ is

$$N_p := \pi(p + \pi(p)) - \pi(p) \geq \pi\left(p + \frac{p}{\log p - 1}\right) - \frac{p}{\log p - 1.1}$$

$$\geq \frac{p + \frac{p}{\log p - 1}}{\log\left(p + \frac{p}{\log p - 1}\right) - 1} - \frac{p}{\log p - 1.1}$$

$$= p \left[ \frac{1 + \frac{1}{\log p - 1}}{\log p + \log\left(1 + \frac{1}{\log p - 1}\right) - 1} - \frac{1}{\log p - 1.1} \right]$$

Since, $1/(\log p - 1) > 0$ for $p > 60184$ and $\log(1 + x) < x$ if $x > 0$,

$$N_p > p \left[ \frac{\frac{\log p}{\log p - 1}}{(\log p - 1) + \frac{1}{\log p - 1}} - \frac{1}{\log p - 1.1} \right]$$

$$= p \left[ \frac{\log p}{(\log p - 1)^2 + 1} - \frac{1}{\log p - 1.1} \right]$$

$$= \frac{(0.9 \log p - 2) p}{\log^3 p - 3.1 \log^2 p + 4.2 \log p - 2.2}$$

$$= \frac{0.9\, p}{\log^2 p - \frac{79}{90} \log p + \frac{911}{405} + \frac{10201}{3645} \frac{1}{\log p - 20/9}}$$

in which $-\frac{79}{90} \log p + \frac{911}{405} + \frac{10201}{3645} \frac{1}{\log p - 20/9} < 0$ for $p > 60184$, thus

$$N_p > \frac{0.9 p}{\log^2 p} > 1.$$

It can be verified that, for $p < 60184$, $\pi(p + \pi(p)) - \pi(p) \geq 1$ (see Column 4 in Table 1).

□

With Theorem 1 and the lower bound of $\frac{x}{\log x - 1.1}$ for the prime counting function by Dusart[11], we have the following corollary:

*Corollary 1: The prime gap of a prime number p is less than $p/(\log p - 1)$ for $p \geq 5$.*

$$g(p) < \frac{p}{\log p - 1.1}$$

*Proof:* From Theorem 1 and Dusart's upper bound[11], we have, for $p > 60184$,

$$g(p) \leq \pi(p) < \frac{p}{\log p - 1.1}.$$

It can be verified that, $g(p) < \frac{p}{\log p - 1.1}$ also holds for $5 \leq p < 60184$. □

Numerical verification (last column in Table 1) indicates that a slightly tighter upper bound

$$g(p) < \frac{p+1}{\log p}$$

holds for all prime numbers. Appendix 1 is the Python code used to obtain the data given in Table 1.

Table 1. (*maximal prime gaps in Column 3 marked by* *)

| $n = \pi(p_n)$ | $p_n$ | $g(p_n)$ | $\pi(p_n+n) - n$ | $(p+1)/\log p$ |
|---|---|---|---|---|
| 1 | 2 | 1* | 1 | 4.3 |
| 2 | 3 | 2* | 1 | 3.6 |
| 3 | 5 | 2 | 1 | 3.7 |
| 4 | 7 | 4* | 1 | 4.1 |
| 5 | 11 | 2 | 1 | 5.0 |
| 6 | 13 | 4 | 2 | 5.5 |
| 7 | 17 | 2 | 2 | 6.4 |
| 8 | 19 | 4 | 1 | 6.8 |
| 9 | 23 | 6* | 2 | 7.7 |
| 10 | 29 | 2 | 2 | 8.9 |
| 11 | 31 | 6 | 2 | 9.3 |
| 12 | 37 | 4 | 3 | 10.5 |
| 13 | 41 | 2 | 3 | 11.3 |
| 14 | 43 | 4 | 2 | 11.7 |
| 15 | 47 | 6 | 3 | 12.5 |
| 16 | 53 | 6 | 3 | 13.6 |
| 17 | 59 | 2 | 4 | 14.7 |
| 18 | 61 | 6 | 4 | 15.1 |
| 19 | 67 | 4 | 4 | 16.2 |
| 20 | 71 | 2 | 4 | 16.9 |
| 21 | 73 | 6 | 3 | 17.2 |
| 22 | 79 | 4 | 4 | 18.3 |
| 23 | 83 | 6 | 4 | 19.0 |
| 24 | 89 | 8* | 6 | 20.1 |
| 25 | 97 | 4 | 5 | 21.4 |
| 26 | 101 | 2 | 5 | 22.1 |
| 27 | 103 | 4 | 4 | 22.4 |
| 28 | 107 | 2 | 4 | 23.1 |
| 29 | 109 | 4 | 4 | 23.4 |
| 30 | 113 | 14* | 4 | 24.1 |
| 99 | 523 | 18* | 15 | 83.7 |
| 154 | 887 | 20* | 21 | 130.8 |
| 189 | 1129 | 22* | 25 | 160.8 |
| 217 | 1327 | 34* | 26 | 184.7 |
| 1183 | 9551 | 36* | 126 | 1042.3 |
| 1831 | 15683 | 44* | 184 | 1623.5 |
| 2225 | 19609 | 52* | 223 | 1984.1 |
| 3385 | 31397 | 72* | 330 | 3032.3 |
| 14357 | 155921 | 86* | 1165 | 13040.1 |
| 30802 | 360653 | 96* | 2386 | 28185.6 |
| 31545 | 370261 | 112* | 2439 | 28877.2 |
| 40933 | 492113 | 114* | 3123 | 37547.4 |
| 103520 | 1349533 | 118* | 7325 | 95608.1 |
| 104071 | 1357201 | 132* | 7349 | 96112.8 |
| 149689 | 2010733 | 148* | 10304 | 138537.5 |
| 325852 | 4652353 | 154* | 21244 | 303028.0 |
| 1094421 | 17051707 | 180* | 65621 | 1024018.3 |
| 1319945 | 20831323 | 210* | 78221 | 1236136.0 |
| 2850174 | 47326693 | 220* | 160910 | 2677972.3 |
| 6957876 | 122164747 | 222* | 373308 | 6560632.0 |
| 10539432 | 189695659 | 234* | 551956 | 9952066.6 |
| 10655462 | 191912783 | 248* | 557801 | 10062250.1 |
| 20684332 | 387096133 | 250* | 1044533 | 19575833.9 |
| 23163298 | 436273009 | 282* | 1163064 | 21930122.7 |

Bertrand's postulate[3] can be proved by Theorem 1.

***Theorem 2: (Bertrand–Chebyshev theorem or Bertrand's postulate):** For $n \geq 1$, the prime gap of the nth prime number $p_n$ is less than the prime number itself, $g(p_n) < p_n$.*

**Proof:** *From Theorem 1, $g(p) \leq \pi(p)$. Let $p_n$ be the nth prime number, we have*

$$g(p_n) \leq \pi(p_n) = n.$$

*Since n is always smaller than the nth prime number, or $n < p_n$, the Bertrand's postulate follows,*

$$g(p_n) \leq n < p_n.$$
□

**References**

[1] Havil, J. Gamma, "Exploring Euler's Constant", p. 170, Princeton, NJ: Princeton University Press, 2003.

[2] H. Cramér, "On the order of magnitude of the differences between consecutive prime numbers," Acta. Arith., 2 (1936) 396-403.

[3]. Joseph Bertrand. Mémoire sur le nombre de valeurs que peut prendre une fonction quand on y permute les lettres qu'elle renferme. Journal de l'Ecole Royale Polytechnique, Cahier 30, Vol. 18 (1845), 123-140.

[4]. P. Tchebychev. Mémoire sur les nombres premiers. Journal de mathématiques pures et appliquées, Sér. 1(1852), 366-390.

[5]. Ramanujan, S. (1919). "A proof of Bertrand's postulate". Journal of the Indian Mathematical Society. 11: 181–182.

[6] Ribenboim, Paulo, The New Book of Prime Number Records, 3rd edition, Springer-Verlag, 1995. New York, NY, pp. xxiv+541, ISBN 0-387-94457-5. p. 252-253

[7]. Oppermann, L. (1882), "Om vor Kundskab om Primtallenes Mængde mellem givne Grændser", Oversigt over det Kongelige Danske Videnskabernes Selskabs Forhandlinger og dets Medlemmers Arbejder: 169–179.

[8] D. Andrica, "Note on a conjecture in prime number theory." Studia Univ. Babes-Bolyai Math. 31 (1986), no. 4, 44--48.

[9] R. C. Baker, G. Harman, and J. Pintz, "The difference between consecutive primes,

II", Proc. London Math. Soc., 83:532562, 2001.

[10] N. J. A. Sloane, R. K. Guy, Sequence A005250 in The On-Line Encyclopedia of Integer Sequences (2020), published electronically at https://oeis.org.

[11] P. Dusart, "Estimates of some functions over primes without RH", arXiv:1002.0442v1 [math.NT] 2 Feb 2010.

Email address：luyaping1@yahoo.com

Appendix 1: Python code

```python
from math import sqrt
def prime_check(Num): #Check a whether a number is a prime or not
    for i in range(2, int(sqrt(Num)) + 1):
        if (Num % i == 0):
            return 0
    return 1

def next_p(p): # Find out the next prime number
    p += 2
    while prime_check(p) == 0:
        p += 2
    return p

def prime_count(Num1,Num2): # Get # of primes in (Num1,Num2] and π(Num2)
    count = 0
    p_end = Num1
    j = Num1 + 2
    for j in range(Num1 + 2, Num2 + 1, 2):
        if prime_check(j) == 1:
            count += 1
            p_end = j
    result = [count, p_end]
    return result

n_max = 100000000
gap_max = 1
ct_last = 1
print(1, 2, gap_max, ct_last)
list1 = [1, 3]
p = 3
p_end_last = 3
for n in range(2, n_max + 1):
    list1 = prime_count(p_end_last, p + n)
    ct = ct_last - 1 + list1[0]
    p_next = next_p(p)
    gap = p_next - p
    if gap > gap_max:
        gap_max = gap
        print(n, p, gap, ct)
    p = p_next
    ct_last = ct
    p_end_last = list1[1]
```